\documentclass[twoside]{article}
\usepackage{amssymb , amsmath , epsfig}
\def\C    {\hbox{\kern 2.7pt\vrule height6pt depth-.35pt \rm \kern-2.7pt C}}

\begin{document}
\title{\vspace{-1in}\parbox{\linewidth}{\footnotesize\noindent
{\sc  }}
 \vspace{\bigskipamount} \\
 Resolvent estimates for 2 dimensional perturbations of plane Couette flow
\thanks{ 
\indent This work was partially supported by
CAPES - BEX1901/99-0 - Bras\'\i lia - Brasil } }
\author{ Pablo Braz e Silva }
\date{}
\maketitle

\begin{abstract}
We present results concerning resolvent estimates for the linear
operator associated with the system of differential equations
governing 2 dimensional perturbations of plane Couette flow. We
prove estimates on the $L_2$ norm of the resolvent of this
operator showing this norm to be proportional to the Reynolds
number $R$ for a region of the unstable half plane. For the
remaining region, we show that the problem can be reduced to
estimating the solution of a homogeneous ordinary differential
equation with non-homogeneous boundary conditions. Numerical
approximations indicate that norm of the resolvent is proportional
to $R$ in the whole region of interest.
\end{abstract}

\newtheorem{theorem}{Theorem}
\newtheorem{lemma}[theorem]{Lemma}
\newtheorem{definition}{Definition}

\renewcommand{\theequation}{\arabic{section}.\arabic{equation}}
\catcode`@=11 \@addtoreset{equation}{section} \catcode`@=12

\section{Introduction}
We consider the Initial Boundary Value problem
\begin{equation} \label{eq1.1}
    \left\{ \begin{array}{l}
        \displaystyle w_t + (w\cdot \nabla )W + (W\cdot \nabla )w + \nabla p = \frac{1}{R}\Delta w + H\\
        \displaystyle \nabla \cdot w = 0 \\
        w(x,0,t) = w(x,1,t) = (0 , 0) \\
        w(x,y,t) = w(x + 1 ,y,t) \\
        w(x,y,0) = (0 , 0)
     \end{array} \right.
\end{equation}
where $w : \mathbb{R}\times [0 , 1] \times [0 , \infty)
\longrightarrow \mathbb{R}^2$ is the unknown function $w(x,y,t) =
(u(x,y,t) , v(x,y,t))$. $W$ is The vector field $W = (y,0)$ and
the Reynolds number $R$ is a positive parameter. The forcing
$H(x,y,t) = (F(x,y,t),G(x,y,t))$ is a given $C^{\infty}$ function
satisfying satisfying $\displaystyle\int_{0}^{\infty} ||H(\cdot ,
\cdot , t)||^2 dt < \infty$ and $\nabla \cdot H = 0 $ for all $(x
, y , t) \in \mathbb{R}\times [0,1] \times [0 , \infty)$. These
equations are the linearization of the equations governing 2
dimensional perturbations $w(x,y,t)$ of $W = (y,0)$, known as
Couette flow, which is a steady solution of
\begin{equation}\label{eq1}
\left\{ \begin{array}{l}
        \displaystyle U_t + (U\cdot \nabla )U +  \nabla P = \frac{1}{R}\Delta U \\
        \displaystyle \nabla \cdot U = 0 \\
        U(x,0,t) = (0 , 0) \\ U(x,1,t)= (1 , 0) \\
        U(x,y,t) = U(x + 1 ,y,t) . \\
     \end{array} \right.
\end{equation}
\noindent The pressure term $p(x,y,t)$ in (\ref{eq1.1}) is
determined in terms of $w$ by the linear elliptic equation
\begin{equation} \label{eq1.2} \left\{\begin{array}{l}
    \displaystyle \Delta p = -\nabla \cdot ((w\cdot \nabla) W) - \nabla \cdot ((W\cdot \nabla) w) \vspace{.1cm}\\
    \displaystyle p_y(x,0,t) = \frac{1}{R}v_{yy}(x,0,t) \vspace{.1cm}\\
    \displaystyle p_y(x,1,t) = \frac{1}{R}v_{yy} (x,1,t).
\end{array} \right.
\end{equation}
We note that $p$ depends linearly on $w$, and is determined up to
a constant. The estimates derived in this paper are independent of
$p$. With $p$ given by the above equation, the solution $w$ of
(\ref{eq1.1}) remains divergence free for all $t\geq 0$. Therefore
we drop the continuity equation and write the problem as
\begin{equation}
\label{eq1.5} \left\{ \begin{array}{l} \displaystyle w_t
  = \frac{1}{R}\Delta w - (w\cdot \nabla )W - (W\cdot \nabla )w - \nabla p+ H =: \mathcal{L}_R w + H \\
   w(x,y,0) = (0,0)
\end{array} \right. \end{equation} \noindent with $w$ satisfying the boundary conditions
described in (\ref{eq1.1}). Note that $\displaystyle
\mathcal{L}_R$ is a linear operator on $L_2 \big( \Omega \big)$,
$\Omega = [0 , 1]\times [0 , 1]$. The $L_2$ inner product and norm
over $\Omega$ are denoted by
$$\begin{array}{lcr}
      \displaystyle \langle w_1 , w_2 \rangle = \int_{\Omega} \overline{w}_1 \cdot w_2 \hspace{.1cm} dxdy
      & ; &
      \displaystyle \|w\|^2 = \langle w , w \rangle .
   \end{array}
$$
\noindent As motivation, consider a general linear evolution
equation
\begin{equation}\label{eq1.3}
\left\{ \begin{array}{l}
   w_t = \mathcal{L} w + H \\
   w(x,y,0) = 0 .
\end{array} \right.
\end{equation}
\noindent Formally, after Laplace transformation with respect to
$t$, the transformed equation is $s \widetilde{w} = \mathcal{L}
\widetilde{w} + \widetilde{H}$ for $\mbox{Re}\,  (s)\geq 0$, where
$\widetilde{w}$ is the Laplace transform of $w$. The solution of
the transformed equation is given by
\begin{equation} \label{eq1.6}
   \widetilde{w} = (s\mathcal{I} - \mathcal{L} )^{-1} \widetilde{H} .
\end{equation}
where $\mathcal{I}$ is the identity operator, and then we
get the solution of (\ref{eq1.3}) by applying the inverse Laplace
transform to $\widetilde{w}$. Moreover, from (\ref{eq1.6}),
 \begin{equation}\label{eq1.7}
   \displaystyle \|\widetilde{w} ( \cdot , s)\|^2 \leq \|(s\mathcal{I} - \mathcal{L} )^{-1}\|^2 \|
   \widetilde{H} (\cdot, s)\|^2 .
\end{equation}
We recall the following definition:
\begin{definition}
Let $\mathcal{L}$ be a linear operator in a Banach space $X$ and
$\cal{I}$ the identity operator. The resolvent set $\rho
(\mathcal{L} ) \subseteq \C$ of $\mathcal{L}$ is the set of all
complex numbers $s$ such that the operator $(s\mathcal{I} -
\mathcal{L})^{-1}$ exists, is bounded and has a dense domain in
$X$. The linear operator $(s\mathcal{I} - \mathcal{L})^{-1}$, for
$s \in \rho (\mathcal{L})$, is said to be the resolvent of
$\mathcal{L}$. The set $\sigma (\mathcal{L}) = \C \setminus \rho
(\mathcal{L})$ is the spectrum of $\mathcal{L}$.
\end{definition}

Therefore, if the complex number $s$ is such that $s \in \rho
(\mathcal{L})$, the formal argument above to express the
transformed function $\widetilde{w}$ in terms of the transformed
forcing $\widetilde{H}$ by (\ref{eq1.6}) is valid and we have the
estimate (\ref{eq1.7}). According to Romanov \cite{R}, the
unstable half plane $\mbox{Re}\,  s \geq 0$ is contained in the
resolvent set of the linear operator $\displaystyle \mathcal{L}_R$
defined in (\ref{eq1.5}) , for all positive Reynolds numbers $R$,
and the eigenvalue of $\mathcal{L}_R$ with largest real part is at
a distance from the imaginary axis proportional to $\displaystyle
\frac{1}{R}$. Our aim in this paper is to estimate the resolvent
constant $\displaystyle \sup_{Re\, s\geq0}\|(s\mathcal{I} -
\mathcal{L}_R)^{-1}\|$. For each $R$, we derive estimates on the
$L_2$ norm of the resolvent for the region $|s|\geq
2\sqrt{2}(1+\sqrt{R})$. For the remaining region, we prove that we
can reduce the problem to estimating the norm of the solution of a
homogeneous ordinary differential equation with non-homogeneous
boundary conditions. This is the main contribution of the paper,
since in principle it simplifies the problem either if one wants
to prove the estimates analytically or use numerical computations.
We perform numerical computations that indicate that the $L_2$
norm of the resolvent is proportional to $R$ also for
$|s|<2\sqrt{2}(1+\sqrt{R})$. We note that, according to Romanov
\cite{R}, this is the best possible dependence of the resolvent on
the Reynolds number. This dependence is better than the 3
dimensional case, since results in Liefvendahl \& Kreiss\cite{L2},
Reddy \& Henningson\cite{RE}, Trefethen {\it et al.}\cite{T},
Kreiss {\it et al.}\cite{K}, also partly based in computations,
indicate the $L_2$ norm of the resolvent to be proportional to
$R^2$. In all the papers above, the computations are performed
using different methods than the one used here.

\section{Resolvent estimates}
\subsection{$s$ bounded away from 0}

\begin{theorem} \label{teorema1}
If $|s| \geq 2\sqrt{2}\left( 1+\sqrt{R}\right)$, then
$\displaystyle \|(s\mathcal{I} - \mathcal{L}_R)^{-1}\|^2 \leq
\frac{8}{|s|^2}(1+\sqrt{R})^2 \leq 1$.
\end{theorem}
\paragraph{Proof:}
After Laplace transformation in $t$, the first equation in
(\ref{eq1.1}) is transformed to

\begin{equation} \label{eq2.1} \begin{array}{l} \displaystyle s \widetilde{w}
  = \frac{1}{R}\Delta \widetilde{w} - (\widetilde{w}\cdot \nabla )W - (W\cdot \nabla )\widetilde{w} - \nabla \tilde{p}+ \widetilde{H} =
  \mathcal{L}_R \widetilde{w} + \widetilde{H} \\
\end{array} \end{equation}
with $\widetilde{w}$ satisfying the boundary conditions in the
space variables described in (\ref{eq1.1}). Taking the inner
product of (\ref{eq2.1}) with $\widetilde{w}$, we have
\begin{equation}  \label{eq2.2}
\displaystyle \langle\widetilde{w},s\widetilde{w}\rangle =
\langle\widetilde{w}, \frac{1}{R}\Delta\widetilde{w}\rangle -
\langle\widetilde{w},(\widetilde{w} \cdot \nabla )W\rangle -
\langle\widetilde{w},(W\cdot \nabla)\widetilde{w} \rangle -
\langle\widetilde{w} , \nabla \widetilde{p}\rangle +
\langle\widetilde{w} ,\widetilde{H}\rangle .
\end{equation}
\noindent Integrating by parts and using the divergence free and
boundary conditions, one can prove that $\langle\widetilde{w} ,
\nabla \widetilde{p}\rangle = 0$ and $\langle\widetilde{w},(W\cdot
\nabla)\widetilde{w} \rangle$ is purely imaginary. Therefore,
\begin{equation}  \label{eq2.3}
s\|\widetilde{w}\|^2 + \frac{1}{R}\|D\widetilde{w} \|^2 =
-\langle\widetilde{w},(\widetilde{w} \cdot \nabla )W\rangle  -
\langle\widetilde{w},(W\cdot \nabla)\widetilde{w} \rangle +
\langle\widetilde{w} , \widetilde{H}\rangle =: P
\end{equation}
where $D\widetilde{w}$ denotes the derivative of $\widetilde{w}$
with respect to the space variables, and $\|D\widetilde{w}\|$ its
Frobenius norm. $P$ satisfies
\begin{equation} \label{eq2.4}
\displaystyle | P | \leq  \|\widetilde{w}\| \|(\widetilde{w} \cdot
\nabla)W\| + \|\widetilde{w}\|\|(W\cdot \nabla )\widetilde{w}\| +
\|\widetilde{w}\| \|\widetilde{H}\|
 \end{equation} \noindent which implies, since $\|D W\| = 1$,
\begin{equation}\label{eq2.5}
  | P | \leq  \|\widetilde{w}\|^2 +
\|\widetilde{w}\|\|D\widetilde{w}\| + \|\widetilde{w}
\|\|\widetilde{H}\| .
\end{equation}
 Since
$\langle\widetilde{w},(W\cdot\nabla)\widetilde{w} \rangle$ is
purely imaginary, taking the real part of (\ref{eq2.3}) gives
\begin{equation}\label{equacao1}
     \displaystyle \mbox{Re}\,  s \|\widetilde{w}\|^2 + \frac{1}{R}\|D\widetilde{w} \|^2 \leq
|\langle\widetilde{w},(\widetilde{w} \cdot \nabla )W\rangle| +
|\langle\widetilde{w} , \widetilde{H}\rangle |
              \leq \|\widetilde{w}\|^2 + \|\widetilde{w} \|\|\widetilde{H}\|
\end{equation}
\noindent and therefore
\begin{equation} \label{eq2.6}
 \displaystyle \mbox{Re}\, s \|\widetilde{w}\|^2 \leq \|\widetilde{w}\|^2 + \|\widetilde{w}
\|\|\widetilde{H}\|.
\end{equation}
Note that (\ref{equacao1}) and (\ref{eq2.6}) are valid for all $s$
satisfying $\mbox{Re}\, s \geq 0$.
To complete the proof, we consider two separate cases:\\

\noindent{\bf Case 1:} $\mbox{Re}\,  s \geq |\mbox{Im}\, s|$. \\ \\
In this case, $\displaystyle \mbox{Re}\,  s \geq
\frac{|s|}{\sqrt{2}}$, and we choose $s$ such that $\displaystyle
\frac{|s|}{2\sqrt{2}} \geq 1$. Using (\ref{eq2.6}),
$$
     \frac{|s|}{2\sqrt{2}} \|\widetilde{w}\|^2  = \left( \frac{|s|}{\sqrt{2}} - \frac{|s|}{2\sqrt{2}}\right) \| \widetilde{w} \|^2 \leq ( \mbox{Re}\, s -
      1 ) \|\widetilde{w}\|^2 \leq
\|\widetilde{w}\|\|\widetilde{H}\|
$$
which implies
\begin{equation}\label{eq2.7}
\displaystyle |s|^2 \|\widetilde{w}\|^2 \leq 8\|\widetilde{H}\|^2.
\end{equation}

\noindent{\bf Case 2:} $|\mbox{Im}\, s| \geq \mbox{Re}\,  s \geq 0$. \\  \\
In this case, $\displaystyle |\mbox{Im}\, s| \geq
\frac{|s|}{\sqrt{2}}$. Taking the imaginary part of equation
(\ref{eq2.3}) gives
\begin{equation} \label{eq2.9}
    \displaystyle |\mbox{Im}\, s| \|\widetilde{w}\|^2 \leq  |P|
\leq \|\widetilde{w}\|^2 + \|\widetilde{w}\|\|D\widetilde{w}\| +
\|\widetilde{w} \|\|\widetilde{H}\| .
\end{equation}
By (\ref{equacao1}),
$$
\displaystyle \frac{1}{R}\|D\widetilde{w}\|^2 \leq
\|\widetilde{w}\|^2 + \|\widetilde{w}\|\|\widetilde{H}\|\leq
\left(\|\widetilde{w}\| + \|\widetilde{H}\|\right)^2
$$
and then
\begin{equation}\label{equacao2}
  \|D\widetilde{w}\| \leq \sqrt{R} \left( \|\widetilde{w}\| +
\|\widetilde{H}\|\right).
\end{equation}
Using this estimate in (\ref{eq2.9}), we get
\begin{equation} \label{equacao3}
\displaystyle |\mbox{Im}\, s| \|\widetilde{w}\|^2 \leq
(1+\sqrt{R})\|\widetilde{w}\|^2 + (1+\sqrt{R})\|\widetilde{w}
\|\|\widetilde{H}\|.
\end{equation}
Choosing $s$ such that $\displaystyle \frac{|s|}{2\sqrt{2}} \geq
1+\sqrt{R}$, we have
$$
     \frac{|s|}{2\sqrt{2}} \|\widetilde{w}\|^2 \leq \left( | \mbox{Im}\, s |-
      (1+\sqrt{R}) \right) \|\widetilde{w}\|^2 \leq (1+\sqrt{R})\|\widetilde{w}\| \|\widetilde{H}\|
$$
and this implies
\begin{equation}\label{eq2.8}
\displaystyle |s|^2 \|\widetilde{w}\|^2 \leq
8(1+\sqrt{R})^2\|\widetilde{H}\|^2 .
\end{equation}
From the two cases, we conclude that
\begin{equation} \label{eq2.12}
|s| \geq 2\sqrt{2}(1+\sqrt{R})  \Rightarrow |s|^2 \| \widetilde{w}
\|^2 \leq 8(1+\sqrt{R})^2\|\widetilde{H}\|^2
\end{equation}
which implies the desired estimate
\begin{equation} \label{eq2.13}
   \|(s\mathcal{I} - \mathcal{L}_R)^{-1}\|^2 \leq
\frac{8}{|s|^2}(1+\sqrt{R})^2 \leq 1,
\end{equation}
valid in the region $|s| \geq 2\sqrt{2}(1+\sqrt{R})$ of the
unstable half plane $\mbox{Re}\, s \geq 0$.

\subsection{$\mathbf{|s| < 2\sqrt{2}(1+\sqrt{R})}$} For this case,
write the problem (\ref{eq1.1}) componentwise:
\begin{equation} \label{eq3.1}
    \left\{ \begin{array}{l}
        \displaystyle u_t + y u_x + v + p_x = \frac{1}{R}\Delta u + F\\
        \displaystyle v_t + y v_x  + p_y = \frac{1}{R}\Delta v + G\\
        \displaystyle u_x + v_y  = 0 \\
        u(x,0,t) = u(x,1,t) = v(x,0,t) = v(x,1,t) = 0  \\
        u(x,y,t) = u(x+1 ,y,t) \\
        v(x,y,t) = v(x+1,y,t) \\
        u(x,y,0) = v(x,y,0) = 0
     \end{array} \right.
\end{equation}
Introduce the stream function $\psi$ by \begin{equation}
\label{eq3.2} \left\{
\begin{array}{l}
     \psi_x = v \\
     \psi_y = -u .
\end{array} \right.
\end{equation}
Laplace transform problem (\ref{eq3.1}) with respect to $t$,
expand in a Fourier series in the periodic direction $x$. The
equations for the transformed functions $\widehat{u}(k,y,s)$,
$\widehat{v}(k,y,s)$, $\widehat{p}(k,y,s)$, $\widehat{F} (k,y,s)$,
$\widehat{G} (k,y,s)$ are
\begin{equation}
\label{eq3.3}
  \left\{ \begin{array}{l}
        \displaystyle s \widehat{u} + i k y  \widehat{u} + \widehat{v} + ik\widehat{p} = - \frac{k^2}{R}\widehat{u} +
        \frac{1}{R}\widehat{u}_{yy}+ \widehat{F}\\
        \displaystyle s \widehat{v} + i k y  \widehat{v} + \widehat{p}_y = - \frac{k^2}{R}\widehat{v} +\frac{1}{R}
        \widehat{v}_{yy}+ \widehat{G}\\
        \displaystyle i k u + \widehat{v}_y  = 0 .
     \end{array} \right.
\end{equation}
Differentiating the first equation in (\ref{eq3.1}) with respect
to $y$, the second with respect to $x$ and subtracting the second
from the first, we eliminate the pressure $p$, and using the
divergence free condition, we see that the transformed stream
function $\widehat{\psi}$ satisfies the following boundary value
problem for an ordinary differential equation with three
parameters $s$, $R$, $k$:
\begin{equation} \label{eq3.4}
  \left\{ \begin{array}{l}
 \displaystyle \frac{1}{R} \widehat{\psi} '''' - \left(s + \frac{2k^2}{R} +iky\right)\widehat{\psi} '' +
 \left(sk^2 +\frac{k^4}{R} +ik^3 y\right)\widehat{\psi} = I \\
  \widehat{\psi} (k,0,s) = \widehat{\psi} (k,1,s) = \widehat{\psi} ' (k,0,s) = \widehat{\psi} ' (k,1,s) = 0
\end{array} \right.
\end{equation}
where $I := \widehat{F}_y - ik\widehat{G}$ and $'$ denotes the
derivative with respect to $y$. We also use the notation
$\displaystyle \mathcal{D} = \frac{\partial}{\partial y}$ , and
write the problem as
\begin{equation} \label{eq3.5}
  \left\{ \begin{array}{l}
  \displaystyle T T_0 \widehat{\psi} = \left( \frac{1}{R} \mathcal{D}^2 - (s + \frac{k^2}{R} + i k y) \right) (\mathcal{D}^2 - k^2)
  \widehat{\psi} = I\\
    \widehat{\psi} (k,0,s) = \widehat{\psi} (k,1,s) = \widehat{\psi}_y (k,0,s) = \widehat{\psi}_y (k,1,s) = 0
\end{array} \right.
\end{equation} \noindent where $\displaystyle T = \frac{1}{R} \mathcal{D}^2 - \left( s +
\frac{k^2}{R} +ik y\right)$ and $\displaystyle T_0 = \mathcal{D}^2
- k^2$ are differential operators depending on the parameters $R$,
$k$ and $s$. To derive the resolvent estimates, we use the
following Lemma:
\begin{lemma} \label{lema1}
If for all $k \in \Bbb Z$ and for all $s \in \C$ such that $|s| <
2\sqrt{2}(1+\sqrt{R})$ the solution $\widehat{\psi} (k,y,s)$ of
(\ref{eq3.4}) satisfies
\begin{equation}\label{equacao}
\left\{ \begin{array}{l}
   \displaystyle k^2 \|\widehat{\psi} (k, \cdot , s)\|^2 \leq
C R^2 (\|\widehat{F}(k, \cdot , s)\|^2 + \|\widehat{G}(k, \cdot , s)\|^2) \\
   \displaystyle \|\widehat{\psi} ' (k , \cdot , s)\|^2 \leq
C R^2 (\|\widehat{F}(k, \cdot , s)\|^2 + \|\widehat{G}(k, \cdot ,
s)\|^2)
\end{array} \right.
\end{equation}
where $C$ is a constant independent of $s$, $R$, $k$,
$\widehat{F}$, $\widehat{G}$, then
$$
\|(s\mathcal{I} - \mathcal{L}_R)^{-1}\|^2 \leq 2C R^2
$$
for $|s| < 2\sqrt{2}(1+\sqrt{R})$.
\end{lemma}
\paragraph{Proof:}
If for all $k \in \Bbb Z$ and for all $s \in \C$ such that $|s| <
2\sqrt{2}(1+\sqrt{R})$
\begin{equation} \label{eq3.6}
   \left\{ \begin{array}{l}
   \displaystyle k^2 \|\widehat{\psi} (k, \cdot , s)\|^2 \leq
C R^2 (\|\widehat{F}(k, \cdot , s)\|^2 + \|\widehat{G}(k, \cdot , s)\|^2) \\
   \displaystyle \|\widehat{\psi} ' (k , \cdot , s)\|^2 \leq
C R^2 (\|\widehat{F}(k, \cdot , s)\|^2 + \|\widehat{G}(k, \cdot ,
s)\|^2)
\end{array} \right.
\end{equation} then by
the definition (\ref{eq3.2}) of the stream function,
\begin{equation} \label{eq3.7}
   \left\{ \begin{array}{l}
    \displaystyle \|\widehat{v} (k , \cdot , s) \|^2  = k^2 \|\widehat{\psi} (k , \cdot , s)\|^2 \leq
C R^2(\|\widehat{F}(k, \cdot , s)\|^2 + \|\widehat{G}(k, \cdot , s)\|^2) \\
    \displaystyle \|\widehat{u} (k , \cdot , s) \|^2 = \|\widehat{\psi} ' (k , \cdot , s)\|^2 \leq
C R^2 (\|\widehat{F}(k, \cdot , s)\|^2 + \|\widehat{G}(k, \cdot ,
s)\|^2)
\end{array} \right.
\end{equation} \noindent and since
$$
\begin{array}{l}
\displaystyle \| \widetilde{u} (\cdot , \cdot , s) \|^2 = \sum_{k=-\infty}^{\infty}\| \widehat{u} (k , \cdot,s)\|^2\\
\displaystyle \| \widetilde{v} (\cdot , \cdot , s) \|^2 = \sum_{k=-\infty}^{\infty} \| \widehat{v} (k ,\cdot,s) \|^2 \\
\displaystyle \|\widetilde{F} ( \cdot , \cdot , s) \|^2 +
\|\widetilde{G} ( \cdot , \cdot , s) \|^2 =\sum_{k =
-\infty}^{\infty}(\|\widehat{F} ( k , \cdot , s) \|^2 +
\|\widehat{G} ( k ,\cdot, s) \|^2) ,
\end{array}
$$
(\ref{eq3.7}) implies
 \begin{equation}
\label{eq3.8}
     \left\{ \begin{array}{l}
    \displaystyle \| \widetilde{u} (\cdot , \cdot , s) \|^2\leq
C R^2 ( \|\widetilde{F} ( \cdot , \cdot , s) \|^2 + \|\widetilde{G} ( \cdot , \cdot , s) \|^2)\\
    \displaystyle \| \widetilde{v} (\cdot , \cdot , s) \|^2
\leq  C R^2 ( \|\widetilde{F} ( \cdot , \cdot , s) \|^2 +
\|\widetilde{G} ( \cdot , \cdot , s) \|^2)
   \end{array}\right.
\end{equation} and this implies $\displaystyle \|(s\mathcal{I} -
\mathcal{L}_R)^{-1}\|^2 \leq 2 C R^2$ for $|s| <
2\sqrt{2}(1+\sqrt{R}).$

  This Lemma shows
that to estimate the norm of the resolvent, it is enough to prove
estimates of the form (\ref{eq3.6}) for the solution
$\widehat{\psi}$ of (\ref{eq3.4}). To prove those estimates for
$\widehat{\psi}$, take the inner product of the differential
equation in (\ref{eq3.4}) with $\widehat{\psi}$ and obtain
\begin{equation} \label{eq3.9}
  \displaystyle \frac{1}{R} \langle \widehat{\psi} , \widehat{\psi} '''' \rangle -
  \langle \widehat{\psi} , (s + \frac{2k^2}{R} +iky)\widehat{\psi} '' \rangle
  + \langle \widehat{\psi} ,  (sk^2 +\frac{k^4}{R} +ik^3 y)\widehat{\psi}\rangle  = \langle \widehat{\psi} , I \rangle .
\end{equation}
Through integration by parts,
\begin{equation} \label{eq3.10}
   \left\{ \begin{array}{l}
 \langle \widehat{\psi} , \widehat{\psi} '''' \rangle = \| \widehat{\psi} '' \|^2 \\
 - \langle \widehat{\psi} , \widehat{\psi} '' \rangle  = \| \widehat{\psi} '\|^2 .
\end{array} \right.
\end{equation}
Therefore equation (\ref{eq3.9}) becomes
$$
  \displaystyle \frac{1}{R}  \| \widehat{\psi} '' \|^2 + \left( \frac{2k^2}{R}+s\right) \|\widehat{\psi} '\|^2 +
  \left( s k^2+\frac{k^4}{R}\right) \| \widehat{\psi} \|^2
  -ik\langle \widehat{\psi} , y\widehat{\psi} '' \rangle
  + ik^3\langle \widehat{\psi} ,   y\widehat{\psi}\rangle  = \langle \widehat{\psi} , I \rangle .
$$
Again through integration by parts, we have $\langle
\widehat{\psi} , y\widehat{\psi} '' \rangle = - \langle
\widehat{\psi} , \widehat{\psi} ' \rangle - \langle
y\widehat{\psi} ' , \widehat{\psi} ' \rangle  $, and since $\psi$
is a real function, $\langle \widehat{\psi} , \widehat{\psi} '
\rangle$ is purely imaginary, $\langle y\widehat{\psi} ' ,
\widehat{\psi} ' \rangle$ and $\langle \widehat{\psi} ,
y\widehat{\psi}\rangle$ are real numbers. Therefore, taking the
real part of the previous equation we get
$$
  \displaystyle \frac{1}{R}  \| \widehat{\psi} '' \|^2 + \left( \frac{2k^2}{R}+\mbox{Re}\, s\right) \|\widehat{\psi} '\|^2 +
  \left( (\mbox{Re}\, s) k^2+\frac{k^4}{R}\right) \| \widehat{\psi} \|^2
   = \mbox{Re} \langle \widehat{\psi} , I \rangle - \mbox{Re} (ik \langle \widehat{\psi} , \widehat{\psi} ' \rangle )
$$
and since $\mbox{Re}\,  s \geq 0$ , this equation implies
\begin{equation} \label{eq3.13}
  \displaystyle \frac{1}{R}  \| \widehat{\psi} '' \|^2 +  \frac{2k^2}{R} \|\widehat{\psi} '\|^2 +  \frac{k^4}{R} \| \widehat{\psi} \|^2
    - |k| |\langle \widehat{\psi} , \widehat{\psi} ' \rangle | \leq | \langle \widehat{\psi} , I \rangle | .
\end{equation}
We separate the analysis into three different cases for $k \in
\Bbb Z$: $\displaystyle |k|
> \sqrt{\frac{R}{\sqrt{2}}}$ , $\displaystyle 0<|k|\leq\sqrt{\frac{R}{\sqrt{2}}}$
and the special case $ k = 0$.

\noindent{\bf Case 1:} $\displaystyle |k| > \sqrt{\frac{R}{\sqrt{2}}}$ \\
For this case, we prove the following:
\begin{theorem}
If $\displaystyle |k| >\sqrt{\frac{R}{\sqrt{2}}}$, then
$$
 \left\{ \begin{array}{l}
   \displaystyle k^2 \|\widehat{\psi} (k, \cdot , s)\|^2 \leq 16 R^2 (\|\widehat{F}(k, \cdot , s)\|^2
   + \|\widehat{G}(k, \cdot , s)\|^2) \\
   \displaystyle \|\widehat{\psi} ' (k , \cdot , s)\|^2 \leq 16 R^2 (\|\widehat{F}(k, \cdot , s)\|^2 +
\|\widehat{G}(k, \cdot , s)\|^2).
\end{array}\right.
$$
\end{theorem}
\paragraph{Proof:}
Since $\displaystyle |\langle \widehat{\psi} , \widehat{\psi} '
\rangle | \leq \|\widehat{\psi}\| \|\widehat{\psi} '\| \leq
\frac{R}{4 |k|} \|\widehat{\psi}\|^2 + \frac{|k|}{ R}
\|\widehat{\psi} '\|^2$, inequality (\ref{eq3.13}) implies
\begin{equation} \label{eq3.14}
   \displaystyle \frac{1}{R}  \| \widehat{\psi} '' \|^2 +  \left( \frac{2k^2}{R}  - \frac{k^2}{R}\right) \|\widehat{\psi} '\|^2 +
 \left( \frac{k^4}{R} - \frac{R}{4}\right) \| \widehat{\psi} \|^2 \leq | \langle \widehat{\psi} , I \rangle |
\end{equation} \noindent and since  $\displaystyle |k| > \sqrt{\frac{R}{\sqrt{2}}}$ , we have $
\displaystyle \frac{k^4}{R} - \frac{R}{4} > \frac{k^4}{2 R}$ and
(\ref{eq3.14}) implies \begin{equation} \label{eq3.15}
   \displaystyle \frac{1}{R}  \| \widehat{\psi} '' \|^2 + \frac{k^2}{R} \|\widehat{\psi} '\|^2 +
  \frac{k^4}{2R} \| \widehat{\psi} \|^2 \leq | \langle \widehat{\psi} , I \rangle | .
\end{equation}
Since equation (\ref{eq3.4}) is linear and the forcing is $I =
\widehat{F}_y - ik\widehat{G}$, it suffices, to get the desired
estimates for the solution $\widehat{\psi}$, to prove estimates
for $\widehat{\psi}_1$ and $\widehat{\psi}_2$, solutions of the
boundary value problems
\begin{equation}
  \left\{ \begin{array}{l}
 \displaystyle \frac{1}{R} \widehat{\psi}_1 '''' - \left(s + \frac{2k^2}{R} +iky\right)\widehat{\psi}_1 '' +
 \left(sk^2 +\frac{k^4}{R} +ik^3 y\right)\widehat{\psi}_1 = \widehat{F}_y \\
  \widehat{\psi}_1 (k,0,s) = \widehat{\psi}_1 (k,1,s) = \widehat{\psi}_1 ' (k,0,s) = \widehat{\psi}_1 ' (k,1,s) = 0
\end{array} \right.
\end{equation} \noindent and
\begin{equation}
  \left\{ \begin{array}{l}
 \displaystyle \frac{1}{R} \widehat{\psi}_2 '''' - \left(s + \frac{2k^2}{R} +iky\right)\widehat{\psi}_2 '' +
 \left(sk^2 +\frac{k^4}{R} +ik^3 y\right)\widehat{\psi}_2 =  ik\widehat{G}\\
  \widehat{\psi}_2 (k,0,s) = \widehat{\psi}_2 (k,1,s) = \widehat{\psi}_2 ' (k,0,s) = \widehat{\psi}_2 ' (k,1,s)
  = 0 ,
\end{array} \right.
\end{equation}
since $\widehat{\psi}$ is given by $\widehat{\psi} = \widehat{\psi}_1 + \widehat{\psi}_2$.\\

\noindent {\bf Estimates for $\widehat{\psi}_1$}: Through
integration by parts,
$$
 | \langle \widehat{\psi}_1 , \widehat{F}_y \rangle | = | \langle \widehat{\psi}_1 ' , \widehat{F} \rangle | .
$$
Using the Cauchy-Schwarz inequality, (\ref{eq3.15}) with forcing
$\widehat{F}_y$ yields
\begin{equation} \label{eq3.16}
\displaystyle \frac{1}{R}  \| \widehat{\psi}_1 '' \|^2 +
\frac{k^2}{R} \|\widehat{\psi}_1 '\|^2 +
  \frac{k^4}{2R} \| \widehat{\psi}_1 \|^2 \leq \|\widehat{\psi}_1 '\|\|\widehat{F}\|
\end{equation}
and since all the terms on the left hand side of this equation are
positive, we have
\begin{equation} \label{eq3.17}
   \displaystyle \frac{k^2}{R}\|\widehat{\psi}_1 '\|^2 \leq \|\widehat{\psi}_1 '\|\|\widehat{F}\|
   \Rightarrow k^4 \|\widehat{\psi}_1 '\|^2 \leq R^2
   \|\widehat{F}\|^2
\end{equation}and
\begin{equation} \label{eq3.18}
    \displaystyle \frac{k^4}{2R} \| \widehat{\psi}_1 \|^2 \leq \|\widehat{\psi}_1 '\|\|\widehat{F}\|
     \Rightarrow
    k^6 \|\widehat{\psi}_1 \|^2 \leq 2 R^2 \|\widehat{F}\|^2
\end{equation} \noindent and \begin{equation} \label{3.19}
  \displaystyle \frac{1}{R}\|\widehat{\psi}_1 '' \|^2 \leq \|\widehat{\psi}_1 '\|\|\widehat{F}\|
  \Rightarrow
  k^2\|\widehat{\psi}_1 '' \|^2 \leq R^2 \|\widehat{F}\|^2 .
\end{equation}
{\bf Estimates for $\widehat{\psi}_2$}: By (\ref{eq3.15}) with
forcing $i k \widehat{G}$ and the Cauchy-Schwarz inequality, we
have
\begin{equation} \label{eq3.20}
\displaystyle \frac{1}{R}  \| \widehat{\psi}_2 '' \|^2 +
\frac{k^2}{R} \|\widehat{\psi}_2 '\|^2 +
  \frac{k^4}{2R} \| \widehat{\psi}_2 \|^2 \leq |k| \|\widehat{\psi}_2 \|\|\widehat{G}\|
\end{equation} \noindent and then, similarly to the argument used to prove the
estimates for $\widehat{\psi}_1$ above, we conclude
\begin{equation} \label{eq3.21}
  \left\{ \begin{array}{l}
  \displaystyle k^6 \|\widehat{\psi}_2 \|^2 \leq 4 R^2 \|\widehat{G}\|^2 \\
      k^4 \|\widehat{\psi}_2 '\|^2 \leq 2 R^2 \|\widehat{G}\|^2 \\
      k^2 \|\widehat{\psi}_2 ''\|^2 \leq 2 R^2 \|\widehat{G}\|^2 .
   \end{array} \right.
\end{equation}
Therefore $\widehat{\psi}$, the solution of problem (\ref{eq3.4}),
satisfies
$$
k^2 \|\widehat{\psi} ''\|^2 + k^4 \|\widehat{\psi} '\|^2 + k^6
\|\widehat{\psi}\|^2
    \leq 16 R^2 (\|\widehat{F}\|^2 + \|\widehat{G}\|^2) = 16 R^2
\|\widehat{H}\|^2 . \hspace{.5cm}  $$

\begin{figure}[!h]
\caption{$\displaystyle \max_{k,\xi}k^2\|\delta_1(k,\cdot
,i\xi)\|^2$ for $1\leq k \leq\sqrt{\frac{R}{\sqrt{2}}}$,
$-2(1+\sqrt{R})\leq \xi\leq 2(1+\sqrt{R})$.} \label{fig1}
\begin{center}
\epsfig{figure = 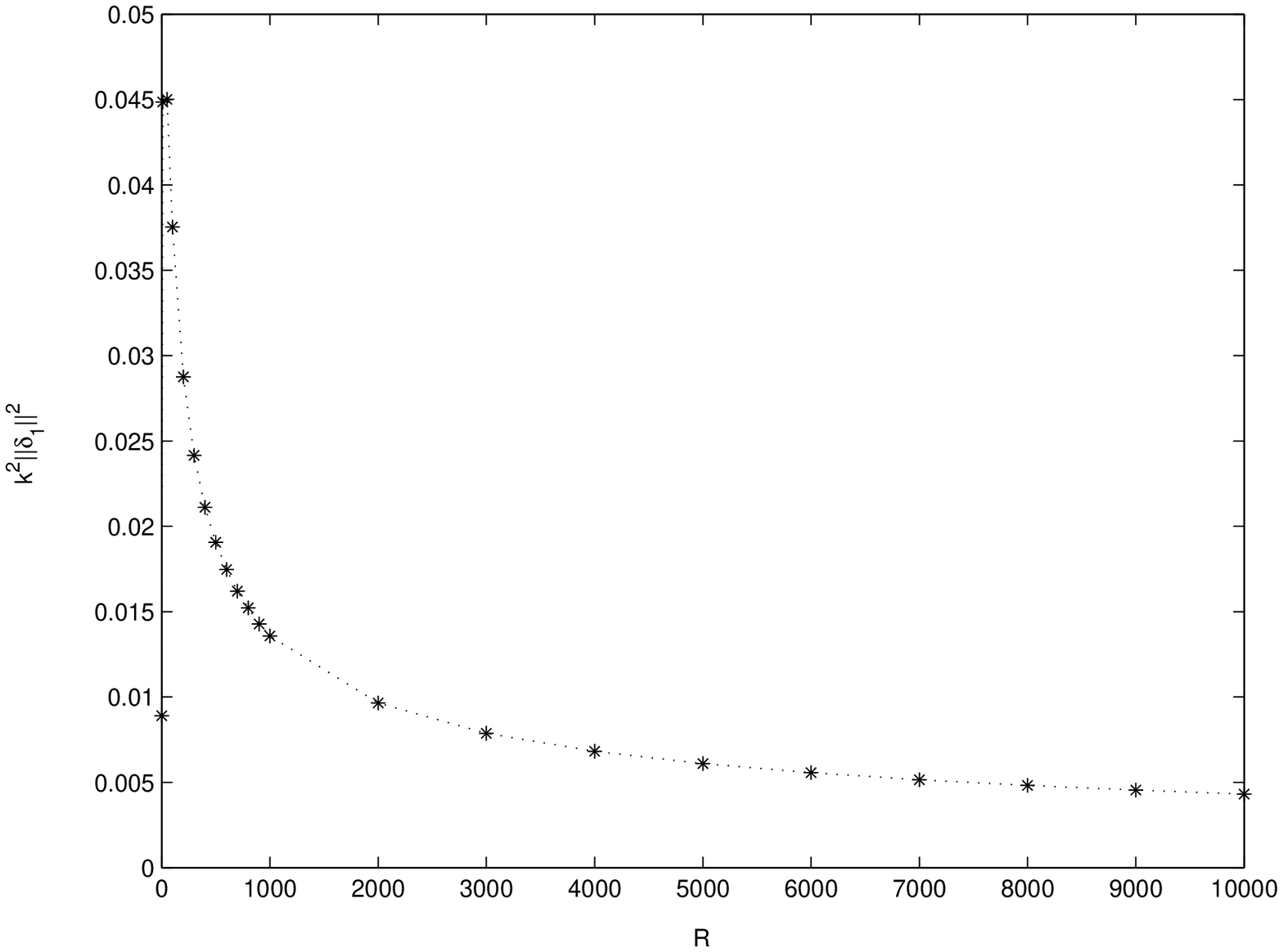,scale=.42}
\end{center}
\end{figure}

\begin{figure}[!h]
\caption{$\displaystyle \max_{k,\xi}\|\delta_1 '(k,\cdot
,i\xi)\|^2$ for $  1\leq k \leq\sqrt{\frac{R}{\sqrt{2}}}$,
 $-2(1+\sqrt{R})\leq \xi\leq 2(1+\sqrt{R})$.}\label{fig2}
\begin{center}
\epsfig{figure = 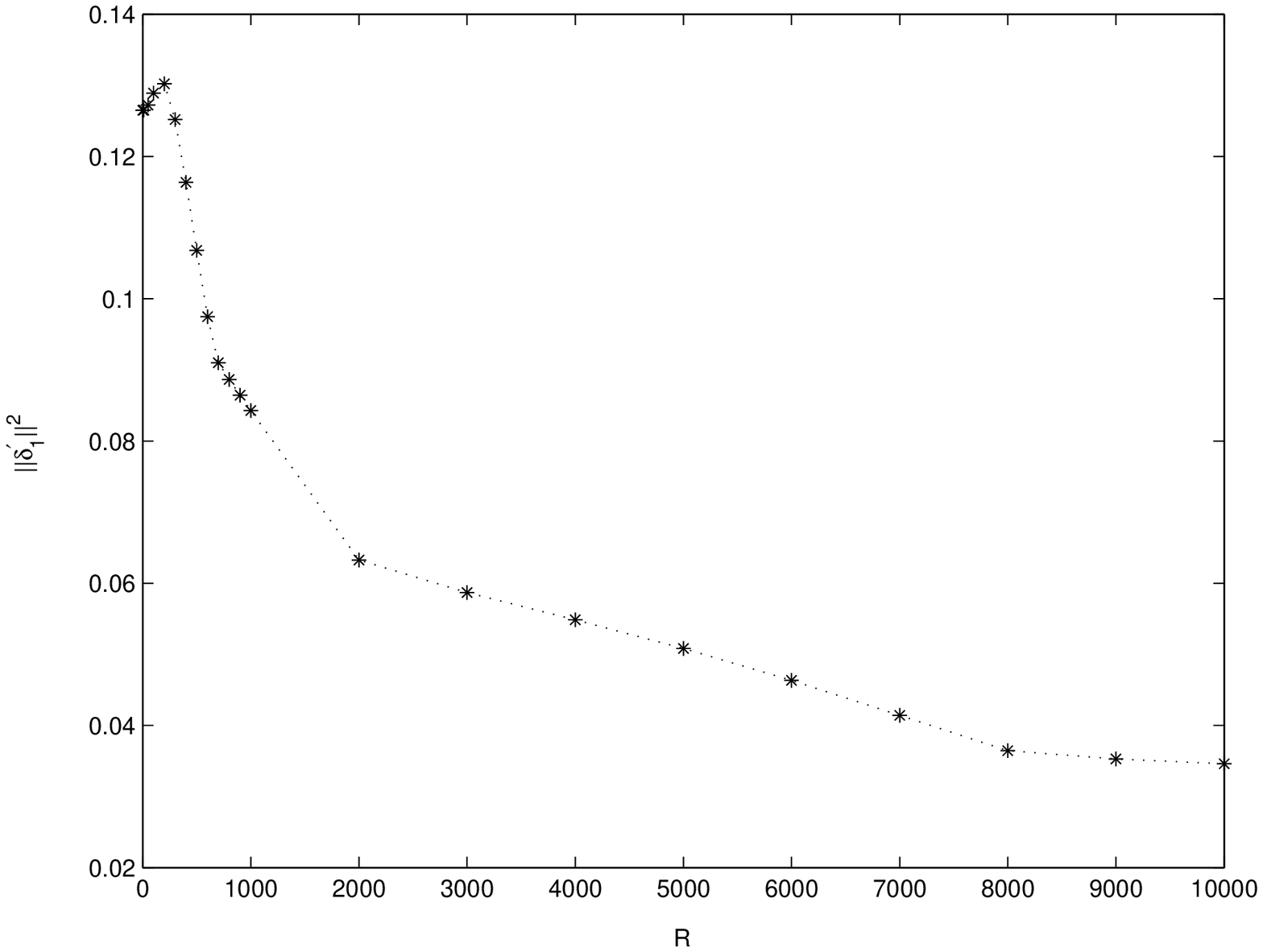,scale=.5}
\end{center}
\end{figure}

\begin{figure}[!h]
\caption{$\displaystyle\max_{k,\xi}k^2\|\delta_2(k,\cdot
,i\xi)\|^2$ for $
 1\leq k \leq\sqrt{\frac{R}{\sqrt{2}}}$,
 $-2(1+\sqrt{R})\leq \xi\leq 2(1+\sqrt{R})$.}\label{fig3}
\begin{center}
\epsfig{figure = 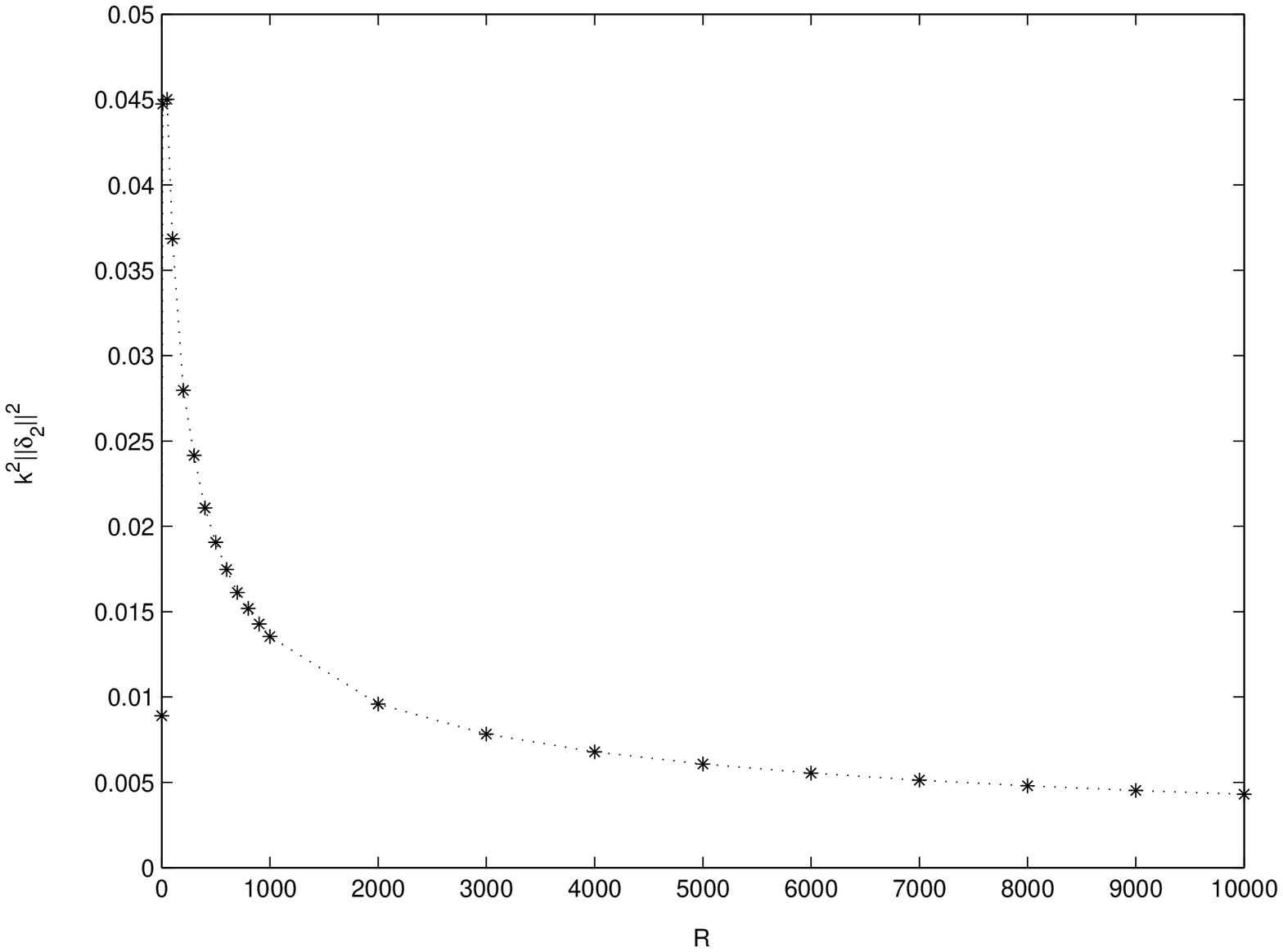,scale=.5}
\end{center}
\end{figure}

\begin{figure}
\caption{$\max\|\delta_2 '(k,\cdot ,i\xi)\|^2$ for $1\leq k
\leq\sqrt{\frac{R}{\sqrt{2}}}$, $-2(1+\sqrt{R})\leq \xi\leq
2(1+\sqrt{R})$.}\label{fig4}
\begin{center}
\epsfig{figure = 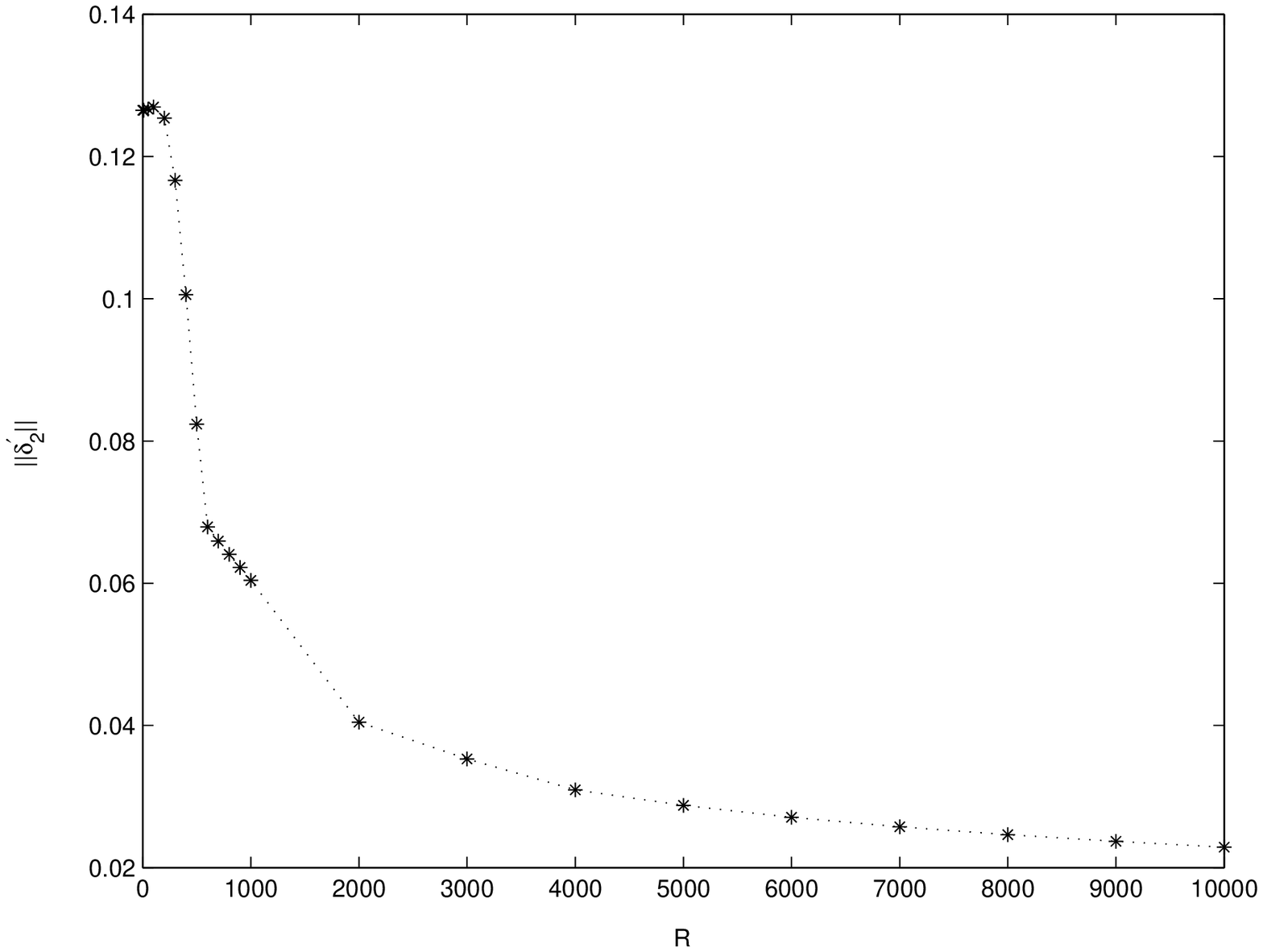,scale=.5}
\end{center}
\end{figure}

\noindent{\bf Case 2:} $\displaystyle 0 < |k| \leq \sqrt{\frac{R}{\sqrt{2}}}$ \\
For this case, we prove that we can reduce the problem to the
study of the solutions of a homogenous 4th order ordinary
differential equation with non-homogenous boundary conditions. We
compute the norms of those functions numerically. The computations
indicate that the norm of the solutions of these simplified
problems do not grow as $R$ grows. This implies the norm of the
resolvent of the operator $\displaystyle\mathcal{L}_R$ to be
proportional to $R$. We begin by noting that we can restrict
ourselves to the case when $s = i \xi$ is purely imaginary. This
is a consequence of the following theorem for holomorphic mappings
in Banach spaces (Chae \cite{C}):
\begin{theorem} \label{the1}
(Maximum Modulus Theorem)
 Let $U$ be a connected open subset of $\C$ and $f : U \longrightarrow E$ a holomorphic
 mapping, where $E$ is a Banach space. If $\|f(z)\|$ has a maximum
 at a point in $U$, then $\|f(z)\|$ is constant on $U$.
\end{theorem}
Applying the theorem to the holomorphic function $f(s) =
(s\mathcal{I} - \mathcal{L}_R)^{-1}$ defined on $\mbox{Re}\, s >
0$, taking values in the Banach space of bounded linear operators
on $L_2(\Omega)$ and noting that (\ref{eq2.13}) implies
$\displaystyle \lim_{|s|\rightarrow\infty}\|f(s)\| = 0$, we
conclude
\begin{equation}
   \displaystyle \sup_{Re s \geq 0} \|s\mathcal{I} - \mathcal{L}_R\| = \sup_{\xi\in\mathbb{R}} \| i\xi\mathcal{I} -
   \mathcal{L}_R\| .
\end{equation}
We note that since we just need to consider $s = i\xi$, $\xi \in
\Bbb R$, and in this case $|s| = |\xi | = |\mbox{Im}\, s|$, the
proof of Theorem \ref{teorema1} is valid for $|s| = |\mbox{Im}\,
s|> 2(1+\sqrt{R})$. This restricts a little bit the parameter
range for the numerical calculations. We now reduce our problem to
the study of solutions of a homogenous ordinary differential
equation. To this end, first consider the second order system
$$
   \begin{array}{lcr}
  \displaystyle T h  = \left( \frac{1}{R} \mathcal{D}^2 - (s + \frac{k^2}{R} + i k y) \right) h = I = \widehat{F}_y - ik\widehat{G}&
  , & h(0) = h(1) = 0 \\
   \displaystyle T_0 g = (\mathcal{D}^2 - k^2) g = h & , & g(0) = g(1) = 0
\end{array}
$$
for $s = i \xi$, $\xi \in \mathbb{R}$. Taking the inner product
with $h$ in the first equation and with $g$ in the second, and
since the boundary conditions for both equations imply that the
boundary terms after integration by parts vanish, we get
\begin{equation}\label{eq3.24}
    \begin{array}{l}
         \|g ''\|^2 + k^2 \| g '\|^2 + k^4 \|g\|^2 \leq C_1\|h\|^2 \\
         \|h'\|^2 + k^2\|h\|^2 \leq C_2 R^2 (\|\widehat{F}\|^2 + \|\widehat{G}\|^2)
    \end{array}
\end{equation}
where $C_1$, $C_2$ are constants independent of $R$, $s$, $k$,
$\widehat{F}$, $\widehat{G}$. Combining those two inequalities, it
follows that
\begin{equation}\label{eq3.25}
   k^2\|g ''\|^2 + k^4 \| g '\|^2 + k^6 \|g\|^2 \leq C R^2 (\|\widehat{F}\|^2 + \|\widehat{G}\|^2)
\end{equation}
where $C$ is a constant independent of $R$, $s$, $k$,
$\widehat{F}$, $\widehat{G}$. Note that $g$ satisfies
\begin{equation}\label{eq3.26}
  \left\{ \begin{array}{l}
     \displaystyle T T_0 g = \left( \frac{1}{R} \mathcal{D}^2 - \left( s + \frac{k^2}{R} + i k y\right) \right) (\mathcal{D}^2 - k^2) g =
            \widehat{F}_y - ik\widehat{G} \\
      g(0) = g(1) = 0 .\\
    \end{array} \right.
\end{equation} \noindent Let $g'(0) = \alpha$ and $g'(1) = \beta$. The numbers
$\alpha$ and $\beta$ can be estimated using a 1-dimensional
Sobolev type inequality. Indeed, $\displaystyle |g'|^2_{\infty}
\leq \|g'\|^2 + \|g'\|\|g''\|$. Using the estimates (\ref{eq3.25})
and that $k \in \Bbb Z$ satisfies $1 \leq |k| $, we conclude that
\begin{equation} \label{eq 3.29}
     \displaystyle |g'|^2_{\infty} \leq C R^2 (\|\widehat{F}\|^2 + \|\widehat{G}\|^2 )
\end{equation} \noindent and therefore \begin{equation} \label{eq3.30}
   \begin{array}{l}
      |\alpha |^2 = |g'(0)|^2 \leq C R^2 (\|\widehat{F}\|^2 + \|\widehat{G}\|^2 ) \\
      |\beta |^2  = |g'(1)|^2 \leq C R^2 (\|\widehat{F}\|^2 + \|\widehat{G}\|^2 )
      ,
   \end{array}
\end{equation}
representing again by $C$ a constant independent of $s$, $R$, $k$,
$\widehat{F}$, $\widehat{G}$. Now, let $\displaystyle \delta (k ,
y , s)$ be the solution of the problem
\begin{equation}\label{eq3.27}
    \left\{ \begin{array}{l}
    \displaystyle T T_0 \delta = \left( \frac{1}{R} \mathcal{D}^2 - \left( s + \frac{k^2}{R} + i k y\right) \right) (\mathcal{D}^2 - k^2) \delta = 0 \\
      \delta (0) = \delta (1) = 0 \\
      \delta ' (0) = \alpha \\
      \delta ' (1) = \beta .
    \end{array} \right.
\end{equation}
\noindent Then, $\displaystyle \widehat{\psi} (k , y , s)  = g (k
, y , s)- \delta (k , y , s) $ is the solution of (\ref{eq3.5}).
Indeed,
$$
  T T_0 \widehat{\psi} = T T_0 (g - \delta) = T T_0 g - T T_0 \delta = \widehat{F}_y - ik\widehat{G}
$$
and
$$ \begin{array}{l}
\displaystyle \widehat{\psi} (k , 0 , s)  = g (k , 0 , s)- \delta (k , 0 , s) = 0 \\
\displaystyle \widehat{\psi} (k , 1 , s)  = g (k ,1,s)- \delta (k,1,s) =0 \\
\displaystyle \widehat{\psi} ' (k , 0 , s)  = g' (k,0,s)- \delta' (k,0,s) = \alpha - \alpha = 0 \\
\displaystyle \widehat{\psi} ' (k , 1 , s)  = g' (k,1,s)- \delta'
(k,1,s) = \beta - \beta = 0 .
\end{array}
$$
Since we already have suitable estimates for $g$ , our problem is
reduced to deriving estimates for the solution of
\begin{equation}\label{eq3.28}
    \left\{ \begin{array}{l}
    \displaystyle T T_0 \delta = \left( \frac{1}{R} \mathcal{D}^2 - \left( s + \frac{k^2}{R} + i k y\right) \right) (\mathcal{D}^2 - k^2) \delta = 0 \\
      \delta (0) = \delta (1) = 0 \\
      \delta ' (0) = \alpha \\
      \delta ' (1) = \beta
    \end{array} \right.
\end{equation}
where $|\alpha |^2 \leq C R^2 (\|\widehat{F}\|^2 +
\|\widehat{G}\|^2 ) $ and $|\beta |^2 \leq C R^2
(\|\widehat{F}\|^2 + \|\widehat{G}\|^2 ) $. Therefore, it remains
to estimate $k^2 \|\delta(k , \cdot , s)\|^2$ and $
\|\delta^\prime (k , \cdot , s)\|^2$ for the parameter range $(k ,
\xi) \in \Bbb Z \times \mathbb{R}$,
$\displaystyle1\leq|k|\leq\sqrt{\frac{R}{\sqrt{2}}}$ , $ |\xi | <
2(1+\sqrt{R})$, where $s = i\xi$. To study this problem
numerically, we simplify it as follows: if $\delta_1$ is the
solution of
\begin{equation}\label{eq3.31}
    \left\{ \begin{array}{l}
    \displaystyle T T_0 \delta_1 = 0 \\
      \delta_1 (0) = \delta_1 (1) = 0 \\
      \delta_1 ' (0) = 1 \\
      \delta_1 ' (1) = 0
    \end{array} \right.
\end{equation}
and $\delta_2$ is the solution of
\begin{equation}\label{eq3.32}
    \left\{ \begin{array}{l}
    \displaystyle T T_0 \delta_2 =  0 \\
      \delta_2 (0) = \delta_2 (1) = 0 \\
      \delta_2 ' (0) = 0 \\
      \delta_2 ' (1) = 1 ,
    \end{array} \right.
\end{equation}
then $\delta = \alpha \delta_1 + \beta \delta_2$ is the solution
of (\ref{eq3.28}). Therefore
$$
k^2 \|\delta(k , \cdot , s)\|^2 \leq 2 |\alpha |^2 k^2\|\delta_1
(k , \cdot , s)\|^2 +2 |\beta |^2 k^2 \|\delta_2 (k , \cdot ,
s)\|^2
$$
and
$$
\|\delta ' (k , \cdot , s)\|^2 \leq 2 |\alpha |^2 \|\delta_1 ' (k
, \cdot , s)\|^2 +2 |\beta |^2  \|\delta_2 ' (k , \cdot , s)\|^2.
$$
Thus we can restrict ourselves to the dependence on $R$ of $k^2 \|
\delta_j (k , \cdot , s)\|^2$ and $\|\delta_j ' (k , \cdot ,
s)\|^2$, for $j =1 , 2$. Moreover, since $\|\psi (k , \cdot ,
s)\|^2 = \| \psi (-k , \cdot , -s)\|^2$, where $\psi$ is the
solution of (\ref{eq3.4}), we can restrict ourselves to $0 \leq
\xi < 2(1+\sqrt{R})$. We performed numerical computations using
the MATLAB 6.0 built-in boundary value problem solver BVP4C, for
the parameter range $1\leq |k|\leq \displaystyle
\sqrt{\frac{R}{\sqrt{2}}}$, $0 \leq \xi \leq 2(1+\sqrt{R})$ and
values of $R$ from $1$ up to $10000$. For $\xi$, we used a mesh
with variable number of points for each $R$. BVP4C makes use of a
collocation method, and we performed computations with different
absolute and relative tolerances, using continuation in the
Reynolds number for the initial guess of the solution. The results
were similar for all cases. Therefore, even though the problem is
stiff for some of the parameter values, the results shown in
figures (\ref{fig1}), (\ref{fig2}), (\ref{fig3}), (\ref{fig4}) in
the following pages should be reliable. For different values of
the Reynolds number $R$, we plot the maximum of $k^2 \| \delta_j
(k , \cdot , s)\|^2$ and $\|\delta_j ' (k , \cdot , s)\|^2$ for
$1\leq k \leq \displaystyle \sqrt{\frac{R}{\sqrt{2}}}$, $0\leq
\xi\leq 2(1+\sqrt{R})$. The results indicate that
\begin{equation} \label{eq3.33}
  \begin{array}{lcr}
  k^2 \| \delta_j (k , \cdot , s)\|^2 \leq 1 &\mbox{and}&
  \|\delta_j ' (k , \cdot , s)\|^2 \leq 1 \hspace{.3cm}\; j=1,2 .
   \end{array}
\end{equation}
Therefore,
$$
   k^2 \|\delta(k , \cdot , s)\|^2 \leq
2 |\alpha |^2 k^2\|\delta_1 (k , \cdot , s)\|^2 +2 |\beta |^2 k^2
\|\delta_2 (k , \cdot , s)\|^2 \leq 4 C R^2 (\|\widehat{F}\|^2 +
\|\widehat{G}\|^2 ),
$$
that is,
\begin{equation}
k^2 \|\delta(k , \cdot , s)\|^2 \leq \widetilde{C} R^2
(\|\widehat{F}(k , \cdot , s)\|^2 + \|\widehat{G}(k , \cdot ,
s)\|^2 )
\end{equation}
and similarly
\begin{equation}
\|\delta ' (k , \cdot , s)\|^2 \leq \widetilde{C} R^2
(\|\widehat{F}(k , \cdot , s)\|^2 + \|\widehat{G}(k , \cdot ,
s)\|^2 )
\end{equation}
 for $\displaystyle
1\leq|k|\leq\sqrt{\frac{R}{\sqrt{2}}}$ , $ |\xi | \leq
2(1+\sqrt{R}) $, where $\widetilde{C}$ is a constant independent
of $R$, $s$, $k$, $\widehat{F}$, $\widehat{G}$. Those inequalities
imply
$$\left\{ \begin{array}{l}
   \displaystyle k^2 \|\widehat{\psi} (k, \cdot , s)\|^2 \leq C R^2 (\|\widehat{F}(k , \cdot , s)\|^2 +
 \|\widehat{G}(k , \cdot , s)\|^2) \\
   \displaystyle \|\widehat{\psi} ' (k , \cdot , s)\|^2 \leq C R^2 (\|\widehat{F}(k , \cdot , s)\|^2 +
 \|\widehat{G}(k , \cdot , s)\|^2) ,
\end{array} \right.
$$
where $C$ is a constant independent of $R$, $s$, $k$,
$\widehat{F}$, $\widehat{G}$. Those are the desired estimates for
$\widehat{\psi}$. \floatsep 1pt
\clearpage
\newpage
\noindent {\bf Case 3:} $\displaystyle k = 0 $\\ For this case,
the differential equation for $\widehat{\psi}$ is reduced to
\begin{equation} \label{eq3.35}
  \left\{ \begin{array}{l}
 \displaystyle \frac{1}{R} \widehat{\psi} '''' - i\xi  \widehat{\psi} ''  = \widehat{F}_y \vspace{.2cm}\\
  \widehat{\psi} (0,0,s) = \widehat{\psi} (0,1,s) = \widehat{\psi} ' (0,0,s) = \widehat{\psi} ' (0,1,s) = 0 .
\end{array} \right.
\end{equation} \noindent By an energy technique, after integration by parts, we get
\begin{equation} \label{eq3.36}
   \displaystyle \frac{1}{R}\|\widehat{\psi} '' \|^2 + i \xi \|\widehat{\psi} '\|^2 = \langle \widehat{\psi} ' , \widehat{F}
\rangle .
\end{equation}
Taking the real part of this equation and using Poincar\'e 's
inequality, we conclude that $\displaystyle \|\widehat{\psi} '
\|^2 \leq  \frac{R^2}{\pi^4}  \|\widehat{F} \|^2$. Applying the
Poincar\'e's inequality again, we get
\begin{equation} \label{eq3.37}
   \displaystyle \|\widehat{\psi}(0, \cdot , s)\|^2 \leq  \frac{R^2}{\pi^6} \|\widehat{F}(0,\cdot ,s)\|^2 .
\end{equation}
which is the desired estimate.

Using Lemma \ref{lema1}, we conclude from the three cases above
that
\begin{equation}
\|(s\mathcal{I} - \mathcal{L}_R)^{-1}\|^2 \leq C R^2
\end{equation}
for $|s| < 2(1+\sqrt{R})$, $C$ a constant independent of $s$, $R$,
$\widehat{F}$, $\widehat{G}$.

\section{Conclusions}
The estimates derived for $s$ bounded away from $0$ and computed
numerically  for the remaining region of the unstable half plane
indicate the $L_2$ norm of the resolvent of the operator
$\mathcal{L}_R$ to be proportional to $R$. Deriving the estimates
analytically for the whole unstable half plane is still an open
problem as far as we know. The method proposed here can be helpful
to solve it. It may also be possible to adapt this method for the
3 dimensional problem, and with this approach it may be possible
to give a satisfactory solution to the problem, or at least to
perform more convincing computations. We hope to address these
questions in the future.

\paragraph{Acknowledgement:} The author wish to thank
Professor Jens Lorenz, from the Department of Mathematics and
Statistics of The University of New Mexico, for suggesting the
problem and for the fruitful discussions about it. He also would
like to thank the anonymous referee, for pointing out a mistake in
a previous proof of Theorem \ref{teorema1}, and for suggestions
that improved the presentation of the paper.

\noindent {\sc Pablo Braz e Silva  } \\
 The University of New Mexico \\
 Department of Mathematics and Statistics \\
 Albuquerque, NM 87131 \\
 e-mail: pablo@math.unm.edu  \smallskip

\end{document}